\documentclass{article}

\usepackage{amsmath,amssymb,amsthm}
\usepackage[margin=1in]{geometry} 
\usepackage{hyperref} 
\usepackage{graphicx} 
\usepackage[english]{babel}
\usepackage{listings}

\title{Definition of the invariant and the relationship with the compounds numbers. Generalisation of the Euler theorem}
\author{Juan Hernández-Toro \\ \texttt{jhtoro@usal.es}}

\date{\today}

\begin{document}
\maketitle 

\begin{abstract}
The purpose of this article is to introduce the concept of invariance and its properties. These properties can be used to check the primality of a number. Combining these properties with the Euler theorem, it is possible to generalize this theorem for all the values of $a^{\varphi(m)}$ where $0 < a < m {\pmod {m}}$ independently if a is co-prime or not with m.
\end{abstract}

\section{Introduction}

The invariant, denoted as I, is a number between 0 and m that accomplishes the following condition: $I^2\equiv I {\pmod {m}}$. This apparently simple property is deeply connected with compound numbers. The principal fact is that, the prime and prime powers don’t have  non-trivial invariants remainders, and for this reason, they can be discriminated. Using this condition, the Euler theorem can be generalized. Applying the invariant properties to the Euler's theorem allows his generalization in the following form:  $a^{\varphi(m)} \equiv I(m) {\pmod {m}}$ independently if a is co-prime or not.

\section{Invariant and anti-invariant}
\subsection{Definition. The invariant}
For each combination of numbers $(a,b)$, $(a',b')$... where  $m=a\cdot b$, $m=a'\cdot b'$...  and $GCD(a,b)=1$, $GCD(a',b')=1$... there are two numbers $I_{1}$ and $I_{2}$ between 0 and m such that $I_{1}^2\equiv I_{1} {\pmod {m}}$ and $I_{2}^2\equiv I_{2} {\pmod {m}}$. These numbers $I_{1}$ and $I_{2}$ are called invariants of $m$.\\

\subsection{Proposition. Each invariants $I_{1}$ and $I_{2}$ are multiple of a or b } \label{Idef1}
In order to accomplish $I^2\equiv I {\pmod {m}}$, $I_{1}$ must be multiple of a and  $I_{2}$ must be multiple of b. 
\subsubsection{Demonstration:}
In order to meet the previous criteria, if  $m=a\cdot b$, the I number can be expressed as $I=\alpha a+x$ and $I=\beta b+y$ Then:
\begin{equation}
 I^{2}=\alpha a \beta b+\alpha ay+\beta bx+xy
\label{Ieq1} 
\end{equation}
and:
\begin{equation}
I^{2}\equiv \alpha ay+\beta bx+xy{\pmod {m}}
\label{Ieq2}
\end{equation}
if $\beta b=\alpha a+x-y$, $\beta\cdot b$ can be substituted in (\ref{Ieq2}) as:
\begin{equation}
I^{2}\equiv \alpha a y+\alpha ax+x^{2} {\pmod {m}} 
\label{Ieq3}  
\end{equation}
There are two possible solutions to achieve 
$I^{2}\equiv I{\pmod {m}}$ :
\begin{itemize}
    \item $x=0$ and $y =1$ that is $I_{1}=\alpha\cdot a$ 
    \item $x=1$ and $y =0$ that is $I_{2}=\beta\cdot b$
\end{itemize}
\subsubsection{Example:}
The number 15 has two combinations of multipliers: a = 1, b = 15, and (a'=3,b'=5). The invariants for (a=1,b=15)  are (1,15), the invariants for (a'=3,b'=5) are (6,10) that is:
\begin{equation}
6^{2}\equiv 6 {\pmod {15}} \text{ and } 10^{2}\equiv 10 {\pmod {15}} 
\end{equation}

\subsection{Proposition. Invariant always exists}\label{Idem1}
If $m=a\cdot b$  and $GCD(a,b)=1$, the invariants $I_{1}$ and $I_{2}$ always exits

\subsubsection{Demonstration:}\
If $m=a\cdot b$  and $GCD(a,b)=1$, must exist two  invariants $I_{1}$ and $I_{2}$. Due to \ref{Idef1},   $I_{1}=\alpha\cdot a$, and $c\equiv 1{\pmod {d}}$. As consequence of the Bezout theorem \cite{Varona} if $GCD(a,b)=1$ then always is possible to find two multiples of a and b that $c\equiv 1{\pmod {d}}$, For $I_{2}$ is possible to follow the same procedure

\subsection{ Definition. The anti-invariant }\label{Idef2}
The Anti-invariant, denoted as A, is the number between 0 and m that accomplishes the following condition  $A^2\equiv m-A{\pmod {m}}$. 

\subsection{Proposition. The existence of an Anti-invariant}
 If a number b accomplish  $b=m-I {\pmod {m}}$, b is an anti-invariant
\subsubsection{Demonstration:}
If $b=m-I$ then: 
\begin{equation}
b^{2}=m^{2}-2\cdot m\cdot b+I^{2}\equiv I=m-b  {\pmod {n}} 
\label{Ieq4}
\end{equation}

\subsubsection{Conclusion:}
If $b=m-I$, then b is an anti-invariant A of m. If $b=m-A$, then b is an invariant I of m.
\subsubsection{Example:}
The number 15 has a total 4 invariants (1,6,10,15). The anti-invariant for each invariant are (14,9,5,0)

 \subsection{Proposition. The existence of anti-invariant, invariant tuples}\label{Idef3}
Each anti-invariant is followed by an invariant. This produces tuples of anti-invariant, invariant consecutive numbers. 
\subsubsection{Demonstration:}

In order to meet the previous criteria, the remainder of (I-1) is calculated.
\begin{equation}
 (I-1)^2=I^{2}-2I+1
\label{Ieq5} 
\end{equation}
but $I^2\equiv I{\pmod {m}}$ for this reason:
\begin{equation}
(I-1)^2\equiv I-2I+1=-I+1\equiv m-(I-1) {\pmod {m}}
\label{Ieq6}
\end{equation}
On the other hand:

\begin{equation}
 (A+1)^2=A^{2}+2A+1
\label{Ieq7} 
\end{equation}
but $A^2\equiv m-A{\pmod {m}}$ then:
\begin{equation}
(A+1)^2 \equiv m-A+2A+1=m+A+1\equiv A+1 {\pmod {m}}
\label{Ieq8}
\end{equation}

\subsubsection{Example:}
The number 15 has four tuples of anti-invariant invariant numbers. (0,1), (5,6), (9,10) and (14,15)

\subsection{Proposition. The existence of trivial anti-invariant, invariant tuple }\label {Idef4}
All the natural numbers m have two invariants and two anti-invariants. These two tuples are easily localized (0,1) and (m-1,m)
\subsubsection{Demonstration:}
This is a particular case of \ref{Idef1}. If $a=1$ and  $m=b$, the two solutions are:
\begin{itemize}
    \item $n=1$ 
    \item $n=m$ 
\end{itemize}
Likewise, using \ref{Idef3} the tuples are defined as (0,1) (m-1,m)
\subsubsection{Conclusion:} 
Concluding this section, all the numbers have 2 invariants, $I=m$ and $I=1$ and 2 anti-invariants $A=m-1$ and $A=0$ (mod m).

\subsection{Definition. The non-trivial anti-invariant, invariant tuple for odd numbers.}\label{Idef5}
If $m=a\cdot b$, m odd, $a\neq 1$, $b\neq 1$ and $GCD(a,b)=1$, then these composite numbers have at least 2 additional anti-invariants, invariants tuples. The localization varies with each number\footnote{The even numbers have additional Trivial invariants}  \footnote{The powers of prime numbers b due to $GCD(a,b)\neq1$ have only 2 anti-invariant, invariant tuples also.}.

\subsection{Proposition. The invariant in superior orders }\label{prop6}
If I is an invariant of m, then $I^s\equiv I{\pmod {m}}$ for all values of s where s is positive and integer.
\subsubsection{Demonstration:}
Suppose that is true for $I^{s-1}$ that is $I^{s-1}\equiv I{\pmod {m}}$ then the remainder of 
\begin{equation}
I^{s}\equiv I^{s-1}\cdot I\equiv I\cdot I\equiv I {\pmod {m}}
 \end{equation}

\clearpage

\section{ Conclusions. The relationship between the compose number, the non-trivial invariant, the Euler theorem, and other hy\-po\-the\-sis and considerations }

\subsection{ The relation between compound numbers and non-trivial invariants}
The existence of a non-trivial invariant automatically indicates that the number is not prime.
\subsubsection{Demonstration:}
The demonstration is almost trivial by \ref{Idef1} and \ref{Idem1}. If it is possible to find an invariant I of m, then by \ref{Idef1} one of the invariant factor, that is $I=\alpha \cdot \beta \cdot \gamma..$, is also factor of m. By \ref{Idem1} if m is a compound number the Invariant always exists. 

\subsection{Using the invariant to test the primality}
The search for invariant could be used to test the primality and even provide additional information related to the factorization. The main advantage is that it is not necessary to check if the number is a power or triangular.\\
Additional to this, The algorithm can be distributed in different machines, can run in both directions and multiple strategies can be used.

\subsection{Algorithm to test the primality using invariants remains}
A very easy algorithm can be created to check the primality and even identifying if is a power of a prime or compound number just searching invariants. This algorithm go from down to up but could be done also in the opposite direction. The procedure is as following\footnote{is not optimize}:
\begin{enumerate}
    \item Create a counter C1. This counter go from 2 to (m-1)/2 and each iteration is increased by 1.
    \item Create a second counter C2 this counter start in 4 and each iteration is increased by 2*C1-1.
    \item If C2 in one iteration is bigger than m then C2=C2-m.
    \item If C2=m then the program stop and m is a power.
    \item If C2=C1 then the program stop and m is a compound number.
    \item Finally if C1=(m-1)/2 an the program doesn´t stop before the program stop and is a prime number.
\end{enumerate}

\subsubsection{Factorization}
The result provide an anti-invariant invariant tuple (c,d). Each of then are multiple of one factor and for this reason factorizing c or d provide one of the factor numbers. Other strategy is multiply $c\cdot d=f$ then $f/m=g$ an the factorization of g give as result all the factors independents of c and d.\\ 
\subsubsection{Code}
Following code implement the previous algorithm. This code is not optimize but give an idea about how easy is.

\begin{lstlisting}[language=Python]


 string = input('please insert a odd number:')
    num = int(string)
    # initialize control variable
    control = int((num - 1) / 2)
    prime = True
    # initialize the other variables
    C1 = 2
    C2 = 4
    # control loop
    while (control >= C1):

        if (C2 > num):
            C2 = C2 - num
        elif (C2 == num):
            print(f'the number is  raised to the a power{C1}:')
            prime = False
            break

        if (C1 == C2):
            print(f'the number is not a prime {C1}:')
            prime = False
            break
        C1 = C1 + 1
        C2 = C2 + 2 * C1 - 1
    if (prime):
        print(f'the number is a prime:')
\end{lstlisting}

\subsection{Generalization of the Euler theorem}\label{Euler}
The Euler theorem can be more general.
Euler's theorem states that, if m and a are co prime positive integers and $\displaystyle \varphi (m)$ is Euler's totient function, then a raised to the power $\displaystyle \varphi (m)$ is congruent to 1 module m.\\
That is:
\begin{equation}
   \displaystyle a^{\varphi (m)}\equiv 1{\pmod {n}}.
\end{equation}\\
For \ref{Idef1} the theorem can be generalized as:\\
If $\displaystyle \varphi (m)$ is Euler's totient function\, 
The number a raised to the power $\displaystyle \varphi (m)$ is congruent to the trivial invariant  module m if m and a are co primes. Then:
\begin{equation}
   \displaystyle a^{\varphi (m)}\equiv I_{trivial}=1{\pmod {m}}.
\end{equation}
If m and a are not co primes then:
\begin{equation}
   \displaystyle a^{\varphi (m)}\equiv I_{not trivial} {\pmod {m}}.
\end{equation}

\subsubsection{Demonstration}
Suppose $m=a\cdot b$, $GCD(a,b)=1$ and the invariant is $I=a^{s}\cdot \beta$ where s is an integer number and $\beta$ is co prime with m and a. Then:
\begin{equation}
   I^{\varphi (m)}= a^{s\cdot \varphi (m)}\cdot \beta^{\varphi (m)} {\pmod {m}}.\label{eqEuler1}
\end{equation}
But $\beta^{\varphi (m)}=1$ by the Euler theorem itself. Substituting:
\begin{equation}
   I^{\varphi (m)}\equiv  a^{s\cdot \varphi (m)} {\pmod {m}}.\label{eqEuler2}
\end{equation}\\

For \ref{prop6} $I^{s\cdot \varphi (m)}\equiv I^{\varphi (m)}\equiv I {\pmod {m}}$. Then  $a^{\varphi (m)}\equiv I {\pmod {m}}$.\\

Is trivial to demonstrate that all the multiples of a  $\delta=a\cdot c$ where c and m are co primes meet
\begin{equation}
   \delta^{\varphi (m)}\equiv I {\pmod {m}}.\label{eqEuler3}
\end{equation}

\subsubsection{Example:}
The Euler´s totien number of 105 is $\varphi (n)=48$ Is easy to check that:
\begin{enumerate}
    \item The multiples of 3, $\beta=3\cdot c$ where c is co-prime with 105 $\beta^{48}\equiv 36 {\pmod {n}}$
    \item The multiples of 5, $\beta=5\cdot c$ where c is co-prime with 105 $\beta^{48}\equiv 85 {\pmod {n}}$
    \item The multiples of 7, $\beta=7\cdot c$ where c is co-prime with 105 $\beta^{48}\equiv 91 {\pmod {n}}$
    \item The multiples of 15, $\beta=15\cdot c$ where c is co-prime with 105 $\beta^{48}\equiv 15\equiv 36\cdot 85 {\pmod {n}}$
    \item The multiples of 21, $\beta=21\cdot c$ where c is co-prime with 105 $\beta^{48}\equiv 21\equiv 36\cdot 91 {\pmod {n}}$
    \item The multiples of 35, $\beta=35\cdot c$ where c is co-prime with 105 $\beta^{48}\equiv 70\equiv 85\cdot 91 {\pmod {n}}$
\end{enumerate}

\subsection{ The invariant, the anti-invariant and the $a^{s}$ cyclic group}\label{cyclic group}

If $m=a\cdot b $ the powers of a $a^{s}$ form a cyclic group. I is equivalent to 1 in the group, and A is equivalent to -1.
\subsubsection{Demonstration:}
If $I_{1}=a\cdot \alpha$ then $I_{1}\cdot a\equiv (A_{1}+1)\cdot a\equiv a$ because the anti invariant is a multiple of b. \\
If $I_{2}=b\cdot \beta$ then $A_{2}\cdot a\equiv (I_{2}-1)\cdot a\equiv -a\equiv m-a$ because the invariant is multiple of b.\\
Also, it is trivial to demonstrate that each $a^s\equiv a\cdot \alpha{\pmod {m}}$ where s is a different rearrangement of the integer numbers between 0 and m-2 and $\alpha$ is an integer number between 1 and m-1.

\subsubsection{Example}
In the following table is possible to find the equivalents of the different subset $\mathbb{Z}_{35_{multiple\; of\;5}}$ multiplication\\
\begin{tabular}{|p{2cm}|p{2cm}|p{2cm}|p{2cm}|p{2cm}|p{2cm}|p{2cm}|}
\hline
 &\textbf{5}&\textbf{10}&\textbf{15}&\textbf{20}&\textbf{25}&\textbf{30} \\ \hline
\textbf{5}&25&15&5&30&20&10\\ \hline
\textbf{10}&15&30&10&25&5&20\\ \hline
\textbf{15}&5&10&15&20&25&30\\ \hline
\textbf{20}&30&25&20&15&10&5\\ \hline
\textbf{25}&20&5&25&10&30&15\\ \hline
\textbf{30}&10&20&30&5&15&25\\ \hline
\end{tabular}
\\
\\
Then the invariant is I=15, the anti invariant is A=20, the inverse of 5 is 10, the inverse of 25 is 30 and the inverse of 15 is 20.

\subsection{The Carmichael number hypothesis }
A hypothesis that could be interesting is the fact that the totien fucntion $\varphi(m)$ is not the minimum function that achieve $\displaystyle a^{\varphi (m)}\equiv I_{trivial}=1{\pmod {m}}$ and can be substituted by $\Omega(m)=lcm((a-1)a^{\alpha-1},(b-1)b^{\beta-1},(c-1)c^{\gamma},...)$ where $m = a^{\alpha} b^{\beta} c^{\lambda} ...$. If the hypothesis is true, then the Carmichael numbers are numbers where p-1 is a multiple of $\Omega(n)$. As example, the $\varphi(m)$ of the first Carmichael number is $\varphi(561) = 320$ because $561 = 3 \cdot 11 \cdot 17$. $560/320$ is not an integer, and consequently, the Euler theorem cannot explain the Carmichael number. On the other side, his, $\Omega$ value is lcm(2, 10, 16) = 80. It is easy to see that 560/80=7, then all the relative primes n $560 \equiv 1 {\pmod {m}}$. In general, all the Carmichael numbers checked accomplish $\frac{p-1}{\Omega}\equiv integer$


\end{document}